\documentclass[10pt]{dimacs-l}

\newcommand{\ext}[2]{#1/#2}
\newcommand{\Gal}[2]{\mathop\mathrm{Gal}(\ext{#1}{#2})}
\newcommand{\kro}[2]{\left( \frac{#1}{#2} \right)}

\newcommand{\Q}{{\bf Q}}

\theoremstyle{definition}
\newtheorem{definition}{Definition}

\begin{document}

\title{Distribution of Values of Real Quadratic Zeta Functions}
\author{Joshua Holden}
\address{Department of Mathematics,
Duke University,
Durham, NC 27708, USA}
\email{holden@math.duke.edu}

\begin{abstract}
    The author has previously extended the theory of regular and
    irregular primes to the setting of arbitrary totally real number
    fields.  It has been conjectured that the Bernoulli numbers, or
    alternatively the values of the Riemann zeta function at odd
    negative integers, are evenly distributed modulo $p$ for every
    $p$.  This is the basis of a well-known heuristic, given by Siegel
    in~\cite{Siegel64},
    for estimating the frequency of irregular primes.
    So far, analyses have shown
    that if $\Q(\sqrt{D})$ is a real quadratic field, then the values
    of the zeta function $\zeta_{D}(1-2m)=\zeta_{\Q(\sqrt{D})}(1-2m)$
    at negative odd integers are also distributed as expected modulo $p$
    for any $p$.  However, it has proven to be very computationally
    intensive to calculate these numbers for large values of $m$.  In
    this paper, we present the alternative of computing
    $\zeta_{D}(1-2m)$ for a fixed value of $D$ and a large number of
    different $m$.
\end{abstract}

\maketitle

\section{Introduction, Conjectures, and Previous Results}

Siegel, in~\cite{Siegel64}, conjectured that the numerators of the
Bernoulli numbers $B_{2m}$ were evenly distributed modulo $p$ for any
odd prime $p$.  This hypothesis was used by Johnson (\cite{Johnson75})
and independently by Wooldridge (\cite{Wooldridge}) to predict the
density of primes with a given index of irregularity, that is such
that $p$ divides a given number of the Bernoulli numbers $B_{2},
\ldots, B_{p-3}$.  It also comes in handy for predicting many other
values that are related to irregular primes, such as the order of
magnitude of the first prime of a given index of irregularity.  (See,
for example, \cite{Wagstaff78}.)

Since $B_{2m}=-\zeta(1-2m)(2m)$, it is equivalent to say that the
values of $\zeta(1-2m)$ are evenly distributed modulo $p$, where
$\zeta(s)$ is the Riemann zeta function.  Of course, the Riemann zeta
function can be generalized to any number field $k$ to get a zeta
function $\zeta_{k}(s)$ associated with the number field $k$.  If $k$
is a totally real number field, there are a number of situations where
the values $\zeta_{k}(s)$ can be seen as analogous to the Bernoulli
numbers, including generalizations of Kummer's Criterion for
predicting when $p$ divides the class number of the $p$-th cyclotomic
field and cases where the equation $x^{p}+y^{p}=z^{p}$ of Fermat's
Last Theorem can be shown not to have solutions in $k$.  (See for
example, \cite{Ernvall79,Ernvall83,Ernvall89}, \cite{HP},
\cite{Greenberg}, \cite{HP}, \cite{Holden98}, and \cite{Kudo}.)

Little or no progress has been made on proving Siegel's hypothesis,
but a great deal of data has been collected, especially in regard to
the prediction of Johnson and Wooldridge.  Specifically, this
prediction says that as $p\to\infty$, the probability that $p$ has
index of irregularity $r$ goes to
$$\left( \frac{1}{2} \right)^{r} \frac{e^{-1/2}}{r!} \enspace .$$
(In addition to the original sources, the details may be found in
Section~5.3 of~\cite{Washington}.)  Note that this prediction does
not rely on the full strength of Siegel's hypothesis, but merely on
the weaker hypothesis that the Bernoulli numbers
are 0 modulo $p$ with probability $1/p$.  The investigations focused
on in this paper relate only to predictions about indices of irregularity
based on this weaker hypothesis.

Wagstaff, in
\cite{Wagstaff78}, computed $u_{r}(x)$, the fraction of primes not
exceeding $x$ with index $r$ of irregularity for each $r$ between 0
and 2 and for all $r\geq 3$ grouped together, and compared this
distribution to the predicted distribution for each multiple $x$ of
1000 up to 125000.  The result of the chi-squared test ``fluctuated
usually between 0.1 and 1.0 and had the value 0.29 at $x=125000$.  It
was 0.03 at $x=8000$'' \cite{Wagstaff78}.  These results correspond to
significance levels of .992, .801, .962, and .999, respectively.
(The significance levels used in this paper correspond roughly to the
probability that the agreement between the observed results and the
predicted results is \emph{not} due to chance.  Statisticians
consider the threshold for considering a result to be not due to
chance to be a significance level of .9 to .95.  Since we are not
actually conducting a valid statistical study in this paper, all of
the statistical results should be taken with a very large grain of
salt.)

Buhler, Crandall, Ernvall, and Mets\"ankyl\"a hold the record for
computations with irregular primes, having found all the irregular
primes below four million as described in \cite{BCEM}.  They do not
seem to have done a chi-squared analysis, but they tabulate the
values of $u_{r}(x)$ for $x=4000000$ and $r$ between 0 and 7.
A chi-squared test using the same methodology as before has the
result 1.02, for a significance level of .796.
Earlier, in \cite{BCS}, Buhler, Crandall, and Sompolski tabulated the
same data for $x=1000000$.  The result of the same chi-squared test is
0.78, for a significance level of .854.

Unfortunately, the
only way to collect data to test Siegel's hypothesis is to
investigate $B_{2m}$ for larger and larger $m$, which is very
computationally intensive.
(See~\cite{Bach}
or~\cite{Fillebrown} for details.)

However, in the more general number field case, there are many more
dimensions to the problem.  We start by restricting our attention to
the case of $k$ an abelian totally real number field.  Then we know that
$$\zeta_{k}(s)=\prod_{\chi\in\hat{G}} L(s, \chi)$$
where $\hat{G}$ is the character group of $G=\Gal{k}{\Q}$ and $L(s,
\chi)$ is the $L$-function associated with the character $\chi$.
Note that $L(s, 1)=\zeta(s)$, so the Riemann zeta function is a factor
of the zeta function for $k$.
(See~\cite{CF}, e.g., for more details.)
Certainly it seems likely that for a fixed (totally real) number field
$k$ and character $\chi$ the values of the numerator of $L(1-2m,
\chi)$ are evenly distributed modulo $p$ as $m$ varies.  (It is known
that these values are rational numbers.)
We also hypothesize that these values for different $\chi$ are
independent, which implies that the numerators of $\zeta_{k}(1-2m)$
are distributed modulo $p$ like the product of $\left| G \right|$
independent integer variables, each of which is evenly distributed
modulo $p$.  We will refer to this as the ``product distribution'',
for lack of a better term.
However, it also is reasonable to conjecture
that for a fixed $m$ the values of $\zeta_{k}(1-2m)$ are distributed
according to the product distribution modulo $p$ as
$k$ varies.  More precisely, if we fix $m$
and the degree of $k$ we expect the values to be distributed
according to the product distribution modulo $p$ as
the discriminant of $k$ varies.  Alternatively, if we fix $m$ and the
discriminant of $k$ we expect the values to be distributed according
to the product distribution modulo $p$ as
the degree varies.

In this paper we will be considering the former situation, with the
degree fixed at $2$, making $k$ a real quadratic field.  We let
$k=\Q(\sqrt{D})$ (where $D$ is the discriminant of $k$) and $\zeta_{D}(s)=\zeta_{k}(s)$.
In this case
$$\zeta_{D}(s)= L(s, 1) L(s, \chi) = \zeta(s) L(s, \chi)$$
where $\chi(s) = \kro{D}{s}$, the Kronecker symbol, where appropriate.

We make the following definitions:

\begin{definition}
Let $k=\Q(\sqrt{D})$ be a real quadratic number field with
discriminant $D$. We
say that an odd prime $p$ is \emph{$k$-regular} (or \emph{$D$-regular}) if $p$ is relatively
prime to $\zeta_{k}(1-2m)$ for all integers $m$ such that $2 \leq  2m
\leq \delta - 2$ and also $p$ is relatively prime to $p
\zeta_{k}(1 - \delta)$, where $\delta = p-1$ unless $D = p$, in which case $\delta
= (p-1)/2$.  The number of such zeta-values that
are divisible by $p$ will be the \emph{index of $k$-irregularity}
(or \emph{index of $D$-irregularity}) of $p$.

Further, we will say that $p$ is \emph{$\chi$-regular} if $p$ is
relatively prime to $L(1-2m, \chi)$ for all integers $m$ such that $2 \leq  2m
\leq \delta - 2$ and also $p$ is relatively prime to $p L(1-2m,
\chi)$.  The number of such $L$-values that
are divisible by $p$ will be the \emph{index of $\chi$-irregularity}
of $p$.
\end{definition}

(See~\cite{Holden98}, \cite{Holden99} or~\cite{mathcomp} for an
explanation of why the definition has exactly this form.)

Saying that the values of $\zeta_{k}(1-2m)$ are distributed according
to the product distribution and that the values of $\zeta(1-2m)$ are
evenly distributed is the same as saying that the values of $L(1-2m,
\chi)$ are evenly distributed modulo $p$.  Then we can make the same
prediction about the indices of $chi$-irregularity that Johnson and
Wooldridge made about the indices of irregularity in the rational
case.  We briefly investigated this issue in~\cite{Holden98}, where
there are tables of the analogue of $u_{r}(x)$ (using the index of
$\chi$-irregularity) for $x=1000$, $r$ from 0 to 4, and $D=5, 8, 12,$
and $13$.  The chi-squared test results are not included, but using
the methodology discussed earlier they are 3.32, 1.74, 1.15, and 2.54.
The corresponding significance levels are .345, .628, .765, and .469,
respectively.  We could total the values of (the analogue of)
$u_{r}(x)$ for the four values of $D$ and compare them to the
predicted values; we might expect that this would give us a better
significance level because of the larger ``sample size''.  However, in
this case the chi-squared result is 3.53 and the significance level is
.316, which is worse than any of the results for the values of $D$
taken separately!  We could instead average the values for the four
values of $D$ and compare them to the predicted values; in this case
the chi-squared result is 0.884 and the significance level is .829,
which is quite good.  However, it is not clear to us what the actual
meaning of these averages are.  Let the reader beware.

The papers
\cite{Holden98} and \cite{mathcomp}
investigate some algorithms, based on formulas given by Siegel in
\cite{Siegel68}, for calculating values of
$\zeta_{D}(1-2m)$.  We show there that while
the best-known algorithms for calculating $\zeta_{D}(1-2m)$ letting $m$
vary take amortized time polynomial in $m$, it is possible to
calculate $\zeta_{D}(1-2m)$ with $m$ fixed and $D$ varying in
amortized time subpolynomial in $D$.  We will use such an algorithm
in the following.

\section{New Results}

In the course of testing the algorithms in~\cite{mathcomp}, we
collected more data in addition to that above. Table~\ref{D5} shows
the number of primes less than 5000 which have $\chi$-index of
irregularity $r$ for various values of $r$ and $D=5$.  We compared the
observed and predicted distributions, using the methodology above,  for
primes below $x$ where $x$ was 1000, 2000, 3000, 4000, and 5000, and
found chi-squared values of 3.32, 5.03, 2.51, 1.73, and 2.10 and
significance levels of .344, .170, .473, .630, and .552, respectively.

\begin{table}[h!]
\caption{Results for $D=5$ and $p<5000$} \label{D5}
\begin{center}
\begin{tabular}{|c|c|c|c|}
\hline
$r$ & number & predicted number & predicted fraction  \\
\hline
0 & 422 & 405.16 &.606531  \\
1 & 186 & 202.58 &.303265  \\
2 & 51 & 50.65 &.075816  \\
3 & 7 & 8.44 &.012636  \\
4 & 2 & 1.06 &.001580  \\
\hline
\end{tabular}
\end{center}
\end{table}

Other data was obtained using the philosophy, described above, of
computing the values of $L(1-2m, \chi)$ for large numbers of $D$ and
relatively small values of $m$.  As in the discussion of $D=5, 8, 12,$ and
$13$ above, we present both the total and the average across the different
discriminants.
Table~\ref{aveD5000} presents the data for all $D<5000$ and $p<100$.  The
chi-squared value for the totals is 81.1 and the significance level
is .000.
The
chi-squared value for the averages is 0.053 and the significance level
is .997.  As before, the actual meaning of these numbers is not clear.

\begin{table}[h!]
\caption{Results for $D<5000$ and $p<100$} \label{aveD5000}
\begin{center}
\begin{tabular}{|c|c|c|c|c|c|}
\hline
$r$ & total & predicted total & average  &  predicted average & predicted  \\
& number & number               & number &  number &            fraction\\
\hline
0 & 21864 & 22068.01 & 14.42 & 14.56 &.606531  \\
1 & 11596 & 11034.01 &7.65 & 7.28 &.303265  \\
2 & 2529 & 2758.50 & 1.67 & 1.82 &.075816  \\
3 &   347 & 459.75 & 0.23 & 0.30 &.012636  \\
4 & 41 & 57.47 & 0.03 & 0.04 &.001580  \\
5 & 7 & 5.75 & 0.005 & 0.004 & .000158 \\
\hline
\end{tabular}
\end{center}
\end{table}

Finally, we computed the values of $L(-1, \chi)$ and $L(-3, \chi)$
for the 303957 discriminants $D$ less than one million, and
calculated the indices of $\chi$-irregularity for the primes 3 and 5.
In this case, since the $p$ involved are very small, we computed the
expected number of primes with given index of irregularity directly
from Siegel's hypothesis,
rather than using the limit as $p \to \infty$.  This takes into
account, for instance, the fact that for $p=3$ the index of
irregularity cannot be more than 1 and for $p=5$ it cannot be more
than 2.  The results are shown in Table~\ref{Dmillion}.  The
chi-squared value for the totals is 24636 and the significance level
is .000.  (In this case the categories used for the chi-squared test
are $r=0$, $r=1$, and $r=2$.) The
chi-squared value for the averages is 0.081 and the significance level
is .960.  Nevertheless, in view of the large amount of data there
seems to be a significant discrepancy between the predicted results
and the observed results, even when the averages are considered.
Coupled with our doubts about the meaning of the averaged numbers,
this leads us to believe that something other than a straightforward
extension of Siegel's hypothesis is influencing the distribution of
the values of $L(1-2m, \chi)$ as the discriminant varies.  We hope to
make the nature of this influence clearer in the future.

\begin{table}[h!]
\caption{Results for $D<1000000$ and $p=3$ or $p=5$} \label{Dmillion}
\begin{center}
\begin{tabular}{|c|c|c|c|c|c|}
\hline
$r$ & total & predicted total & average  &  predicted average & predicted  \\
& number & number               & number &  number &            fraction\\
\hline
0 & 338966 & 397170.48 & 1.115177 & 1.306667 &0.653333 \\
1 & 252832 & 198585.24 & 0.831802 & 0.653333 &0.326667 \\
2 & 16116 & 12158.28 &   0.053021 & 0.040000 &0.020000  \\
\hline
\end{tabular}
\end{center}
\end{table}

\section{Other Conjectures and Future Work}

A number of other conjectures about irregular primes could be extended
to the setting which we have presented. One such,
mentioned by Wagstaff in~\cite{Wagstaff78} without
attribution, is that the irregular primes are evenly distributed
across the possible residue classes modulo $n$ for every positive
integer $n$.  To be exact, one expects that the ratio of irregular
primes in a residue class to all odd primes in the class to be the
same for each possible residue class.
Wagstaff investigated this for primes below 125000 and
$3\leq n \leq 37$ (and also $n=59$ and $67$) and gave tables for $n=3, 4,$
and $5$.  He does not give the chi-squared results, but they are
0.107, 0.060, and 2.420, respectively, with significance levels of
.744, .806, and .490.
We conjecture that the  $\chi$-irregular primes are also distributed
evenly, and intend to investigate this using the data we have collected.

Wagstaff also mentions a conjecture due to Wooldridge about the
distribution of the numbers $2k/p$ for which $(p, 2k)$ is an
irregular pair; i.e. $p$ divides $B_{2k}$.  Wooldridge conjectures
that these numbers have a uniform distribution in the interval
$(0,1)$.  We also intend to investigate whether these numbers are
uniformly distributed for $\chi$-irregular pairs.

Two other conjectures about irregular primes which are widely believed
cannot really be checked using statistical methods.  It is thought
that there are primes with arbitrarily large indices of irregularity;
the largest observed so far is 7, for $p=3238481$ (see~\cite{BCEM}).
The largest index of $\chi$-irregularity yet observed is $5$, which
occurred $7$ times among the primes less than 100 using discriminants
less than 5000.  It is also likely that there are Bernoulli numbers
divisible by arbitrarily large powers of $p$.  Although as yet
$p^{2}$ has not been seen to divide $B_{2m}$, as Washington says,
``there does not seem to be any reason to believe this in
general''~\cite{Washington}.
On the other hand, powers as large as
$3^{7}$ have been observed to divide $L(1-2m, \chi)$; this happens at
$D=3869$ and again at $D=3937$.

Studies such as~\cite{BCEM} also frequently investigate applications
such as Fermat's Last Theorem, Vandiver's conjecture, and Iwasawa's
cyclotomic invariants.  All of these could be modified to the
situation explored here.  Additional data would have to be collected
in order to study these applications; we are not yet certain whether
new algorithms would have to be implemented.

\newcommand{\SortNoop}[1]{}

\end{document}